\documentclass[12pt]{article}
\usepackage{amsmath}
\usepackage{amssymb}

\newcommand{\C}{\mathbb C}
\newcommand{\Z}{\mathbb Z}
\def\N{\mathbb N}

\def\endproof{\hfill{\vrule height4pt width6pt depth2pt}  \bigskip }

\def \proclaim#1{\smallskip\noindent{\bf\ignorespaces#1\unskip.}%
 \bgroup\it\space\ignorespaces}
\def \endproclaim{\par\egroup\smallskip}

\begin{document}

\title{\bf Polarization of an inequality  }
\date{}
\author{Ivo Kleme\v{s} }
\maketitle
{\centerline
{\it
\noindent Department of Mathematics and Statistics,
\noindent 805 Sherbrooke Street West,} }

{\centerline
{\it
\noindent McGill University,
\noindent Montr\'eal, Qu\'ebec,
\noindent H3A 2K6,
\noindent Canada.
}}
\medskip

{\centerline{  Email:  klemes@math.mcgill.ca }}

\bigskip

\bigskip


\noindent {\it Abstract.}
We generalize a previous inequality related to a sharp version of
the Littlewood conjecture on the minimal
$L_1$-norm of $N$-term exponential sums $f$ on the unit circle.
The new result concerns replacing the expression $\log(1+t|f|^2)$
with $\log \left(\sum_{k=1}^K \ t_k|f_{k}|^2 \right)$.
The proof occurs on the level of finite Toeplitz matrices,
where it reduces to an inequality between their polarized determinants
(or ``mixed discriminants").

\vfill
\noindent {\it Acknowledgement: } A part
of this research was carried out while the author was on leave at the Department 
of Mathematics and Statistics, University of Victoria, Canada.
The author is grateful to the department for its hospitality.

\noindent {\it A.M.S. Mathematics Subject Classifications:} 42A32 (15A42).

\bigskip

\noindent {\it Date:} 8 May 2009. {\it Updated:} 6 January 2011.

\eject

\small
\centerline{\bf \S 1. Introduction}
\bigskip

We will prove a simple generalization of
the main results of \cite{k1}.
 For each integer $N \geq 1$,
let $\mathcal{L}(N)$ denote the collection of all complex polynomials $f$
of the form $f(z) = c_0 + c_1 z + c_2 z^{2} + \dots +
c_{N-1} z^{N-1}$ where each $c_j \in \{1,-1\}$ (``Littlewood polynomials" ; see
\cite{bor}). Let $\mathcal{ \widetilde{L}}(N)$ denote the collection of all complex polynomials $f$
of the form $f(z) = c_0 + c_1 z^{m_1} + c_2 z^{m_2} + \dots +
c_{N-1} z^{m_{N-1}}$ where each  coefficient $c_j$ is a complex number with $|c_j| \geq 1$,
and $0<m_1< m_2 < \dots < m_{N-1}$ are integers. Clearly $\mathcal{{L}}(N)
\subset \mathcal{ \widetilde{L}}(N)$.
 Define $D_N \in \mathcal{L}(N)$ by $D_N(z) = 1 +  z + z^2 + \dots +  z^{N-1}$.
 Define the $1$-norm
$||f||_1$ on the unit circle by $ ||f||_1 := \int_0^{2\pi} |f(e^{i\theta})| d\theta/2\pi$.
The Littlewood conjecture concerning $\mathcal{ \widetilde{L}}(N)$ was
that there is an absolute constant $C>0$
such that for all $N$ and all $f \in \mathcal{ \widetilde{L}}(N)$, $ ||f||_1 \geq C||D_N||_1$,
and was proved in \cite{K81} and \cite{MPS}, independently. The ``sharp" Littlewood
conjecture is that one can take $C=1$, and this remains open. The main result of
\cite[Theorem 1.2]{k1} concerns only the smaller family $\mathcal{{L}}(N)$, and states that for all $N \in \N,
f \in \mathcal{{L}}(N)$ and $t>0$,
\begin{equation}
\label{log0}
 \int_0^{2\pi} \log \left(1 + t|D_N(e^{i\theta})|^2 \right) d\theta  \leq
\int_0^{2\pi} \log \left(1 + t|f(e^{i\theta})|^2 \right) d\theta .
\end{equation}
As discussed in \cite{k1},
this implies the sharp Littlewood conjecture  for $f \in \mathcal{{L}}(N)$,
as well as similar sharp $p$-norm inequalities (see below) for the range $0 < p \leq 4$,
by means of some simple integrations over the $t>0$.
The result in the present paper still concerns only the smaller family $\mathcal{{L}}(N)$,
but we generalize (\ref{log0}) in the sense of giving a ``vectorized" or ``polarized" version
of it, as follows:

\bigskip

\proclaim {1.1 Theorem } Let $K \geq 1$ and let $N_1, N_2, \dots , N_K \geq 1$
be given integers. Then for any $f_k \in \mathcal{L}(N_k)$ and any real $t_k > 0, \
1 \leq k \leq K$, we have the inequality
\begin{equation}
\label{log}
\int_0^{2\pi} \log \left(\sum_{k=1}^K \ t_k|D_{N_k}(e^{i\theta})|^2 \right)
 d\theta  \quad \leq \quad
 \int_0^{2\pi} \log \left(\sum_{k=1}^K \ t_k|f_{k}(e^{i\theta})|^2 \right)
 d\theta \ .
\end{equation}
\endproclaim
\bigskip

\noindent The special case $K=2, \ N_1 = 1$ is the old result (\ref{log0}).
As in \cite{k1}, the above theorem immediately implies some $p$-norm inequalities
in $L_p(d\theta)$ in the range $0 < p \leq 4$:
\begin{equation}
\label{L02}
\left|\left| \left(
\sum_{k=1}^K \ t_k|D_{N_k}(e^{i\theta})|^2
\right)^{\frac{1}{2}}\right|\right|_{p} \leq
\left|\left| \left(
\sum_{k=1}^K \ t_k|f_{k}(e^{i\theta})|^2
\right)^{\frac{1}{2}}\right|\right|_{p},  \quad 0 < p \leq 2,
\end{equation}
\begin{equation}
\label{L24}
\left|\left| \left(
\sum_{k=1}^K \ t_k|D_{N_k}(e^{i\theta})|^2
\right)^{\frac{1}{2}}\right|\right|_{p} \geq
\left|\left| \left(
\sum_{k=1}^K \ t_k|f_{k}(e^{i\theta})|^2
\right)^{\frac{1}{2}}\right|\right|_{p},  \quad 2 \leq p \leq 4.
\end{equation}
(To deduce this, we replace $K $ by $K+1$ in the theorem,
put $N_{K+1}=1$, and then use certain integral identities for
the power functions $x^p$ in terms of $\log(1+tx),   0<t<\infty$,
as in \cite[p. 211-212, Lemma 11.1, Ch. 4]{GK}, \cite[Theorem 4]{MS},
 or \cite[page 9]{k1}.)

 The proof of Theorem 1.1 is essentially the same as the proof of
 the special case (\ref{log0}) in \cite{k1}, as will be seen in \S 2.
 The basic lemma
 is again the total unimodularity of $(0,1)$ ``interval matrices'' $M$
(whose intervals of 1's occur in their columns, for instance). The
 only new step is to invoke this fact for the general case of a ``polarized determinant"
$D_n(A_{k_1}, A_{k_2}, \dots, A_{k_n})$ of several $n \times n$
 Gram matrices $A_k = M_k{M_k}^*$, instead
 of only for cases of the type $D_n(I, I, \dots,I, A,A, \dots, A )$ with only the two matrices
$I$ and $A=MM^*$, as was implicitly done in \cite{k1}.
The ``polarized determinant" $D_n(A_{1}, A_{2}, \dots, A_{n})$ of the $n \times n$
matrices $A_i$ can be defined as
$\frac{1}{n!}$ times the
coefficient of $ x_1 \dots x_n$ in
$\det(x_1A_1+\dots + x_nA_n),$ where the $x_i $
are scalars. It has traditionally been called the ``mixed discriminant''
and has been useful in work on the van der Waerden conjecture
\cite{bap}, \cite{gur}. It was also used in \cite{k2} in connection with certain matrices $M$
having complex entries of modulus $\geq 1$, or 0. In this paper we
implicitly use the idea of polarized determinant, but we omit
explicit use of the notation $D_n(A_{1}, A_{2}, \dots, A_{n})$ in the
formal lemmas and proofs.
\bigskip

\centerline{\bf \S 2. Proof of Theorem 1.1. }
\medskip

\proclaim {2.1 Lemma } Let $K, n \in \N$, and for $k=1,\dots, K$ let $M_k$ be
any (rectangular) $n \times m_k$
matrices over $\C$. If $x_k$ are scalars and $A_k := M_k ({M_k}^*) $,  then
\begin{equation}
\label{det1}
\det( x_1A_1 + \dots + x_K A_K) =
 \sum_{n_1+\dots +n_K = n}  \gamma(n_1, \dots, n_K)
\ x_1^{n_1} \dots x_K^{n_K},
\end{equation}
where the coefficients $\gamma(n_1, \dots, n_K)$ are given by
\begin{equation}
\label{det2} \gamma(n_1, \dots, n_K) \ = \
 \sum_{(S_1,\dots, S_K)} \ \left|  \det \big(\ S_1 \ \big| \
 \dots \ \big| \ S_K\ \big)\ \right|^2
\end{equation}
where each $S_k$ denotes an $n \times n_k$ matrix obtained by choosing some
$n_k$ columns of $M_k$ (that is, from $n_k$ distinct column indices) ,
$\big(\ S_1 \ \big| \ \dots \ \big| \ S_K\ \big)$ denotes the
$n \times n$ matrix consisting of
the $K$ blocks $S_1, \dots , S_K$ ,
and the sum is over all ordered $K$-tuples $(S_1,\dots, S_K)$ of such choices
(an empty sum being zero by convention).
\endproclaim
\medskip

\noindent {\bf Proof.} This is a known result \cite{bap}. To prove it,
 consider the two ``block'' matrices $B$ and $C$ defined by $B:= \big(\ {x_1} M_1 \ \big| \
 \dots \ \big| \ {x_K} M_K \ \big)$ and $C:= \big(\  M_1 \ \big| \
 \dots \ \big| \  M_K \ \big)$
Clearly, $BC^*=
x_1A_1 + \dots + x_K A_K$. Now apply the Binet-Cauchy theorem
to expand $\det ( BC^* )$.
\endproof

\noindent
{\bf Remark 1.} In the notation of polarized determinants $D_n$, the
above coefficients are given by
\begin{equation}
\label{gam}
\gamma(n_1, \dots, n_K)=
\frac{n!}{n_1! \dots n_K!} D_n(A_1^{\langle n_1 \rangle}, \dots,
A_K^{\langle n_K \rangle}),
\end{equation}
where $A_1^{\langle n_1 \rangle}$ means $A_1, \dots , A_1$ repeated $n_1$
times, etc.  \cite{bap}.

We review the following facts already used in \cite{k1}:
\medskip

\proclaim {2.2 Lemma \cite{FG}} Let $S$ be a square matrix with
entries in $\{0,1\}$ such that in each column the $1$'s occur in consecutive
row positions (i.e. in an ``interval"). Then:  {\bf (i)} $\det S \in \{-1,0,1\}$.
{\bf (ii)} If $S'$ is a matrix with integer entries
satisfying $S' = S$ {\rm mod 2}, then $ |\det S'| \geq |\det S| $.
\endproclaim
\noindent {\bf Proof.} For (i) see \cite[p. 853]{FG} or \cite[Lemma 2.3]{k1}.
(ii) follows from (i) by noting that $ \det S' = \det S $ mod 2.
\endproof

\proclaim {2.3 Corollary } Let $K, n \in \N$, and for $k=1,\dots, K$
 let $M_k$ be $n \times m_k$
matrices with
entries in $\{0,1\}$ such that in each column the $1$'s occur in consecutive
row positions. Let $M_k'$ be $n \times m_k$
matrices with integer entries such that $M_k' = M_k$ {\rm mod 2} for each $k$.
If $t_k \geq 0$ are scalars and $A_k := M_k ({M_k}^*),\ A_k' := M_k' ({M_k'}^*) \ $, then
\begin{equation}
\label{det3}
\det( t_1A_1 + \dots + t_K A_K) \leq \det( t_1A_1' + \dots + t_K A_K') \ .
\end{equation}

\endproclaim
\noindent {\bf Proof.} Expand both sides of (\ref{det3}) using Lemma 2.1.
Then we have
\begin{equation}
\label{det4}
 \left| \det \big(\ S_1 \ \big| \
 \dots \ \big| \ S_K\ \big)\ \right|^2 \leq \left| \det \big(\ S_1' \ \big| \
 \dots \ \big| \ S_K'\ \big)\ \right|^2
 \end{equation}
  for each of the corresponding
 terms in (\ref{det2}), by Lemma 2.2 applied to the matrices
 $S:= \big(\ S_1 \ \big| \ \dots \ \big| \ S_K\ \big)$ and
$S':= \big(\ S_1' \ \big| \ \dots \ \big| \ S_K'\ \big)$  (which correspond
to the same $K$-tuple of choices of column indices).
\endproof

\noindent
{\bf Remark 2.} For an $n \times n$ Hermitian matrix $A \geq 0$ with eigenvalues $\lambda_i$,
let $||A||_p$ denote the $l_p$-norm $(\lambda_1^p + \dots + \lambda_n^p)^{1/p}$.
Corollary 2.3 implies $l_p$-norm inequalities
for the two matrices on either side of (\ref{det3}),
in the same manner as discussed after Theorem 1.1 (by replacing $K$ with
$K+1$ and taking $M_{K+1} =M_{K+1}'= I $, the $n \times n$ identity matrix):
\begin{equation}
\label{p1}
 \big| \big|t_1A_1 + \dots + t_K A_K\big| \big|_p  \leq
 \big| \big|t_1A_1' + \dots + t_K A_K' \big| \big|_p \ ,  \quad 0 < p \leq 1,
 \end{equation}
 and, if each $M_k'$ has its entries specifically in $\{-1,0,1\}$, then also
 \begin{equation}
\label{p2}
\big| \big|t_1A_1 + \dots + t_K A_K\big| \big|_p  \geq
 \big| \big|t_1A_1' + \dots + t_K A_K' \big| \big|_p \ ,  \quad 1 \leq p \leq 2.
 \end{equation}
The extra assumption is needed for (\ref{p2}) since that
part of the implication relies upon the $l_1$-norms of both sides being the
same.
\medskip

\noindent {\bf Proof of Theorem 1.1.} The proof is the same as the one in
\cite[Theorem 1.2]{k1},
except that here we use Corollary 2.3 in one of the steps. We will repeat the details
for completeness: Let $ \psi(\theta) =
\sum_{k=1}^K \ t_k|f_{k}(e^{i\theta})|^2$,
and for each $n \in \N$ let $T(n,\psi)$ be the $n \times n$ Toeplitz matrix
$T(n,\psi)_{ij} = \widehat{\psi}(j-i), \ 1 \leq i,j \leq n$,
where $\widehat{\psi}(m)$ is the usual Fourier coefficient, $\widehat{\psi}(m)=
\int_0^{2\pi} \psi(\theta) e^{-im\theta} d\theta/2\pi$.
By a theorem of Szeg\"o \cite[\S 5.1, pp. 64-65]{GS},
$(\det T(n,\psi) )^{1/n} \to  \exp(\int_0^{2\pi} \log \psi(\theta) d\theta/2\pi)$
as $n \to \infty$. Fix $n$. It is easy to check that $T(n,\psi) = t_1A_1' + \dots + t_K A_K'$
where each $A_k' = M_k' ({M_k'}^*)$ with $M_k'$ being the $n \times (n+N_k-1)$
Toeplitz matrix $(M_k')_{ij} = \widehat{f_k}(j-i)$, where $\widehat{f_k}(m)=
\int_0^{2\pi} f_{k}(e^{i\theta}) e^{-im\theta} d\theta/2\pi,  \ m \in \Z$; the coefficient
of $z^m$ in the polynomial $f_k(z)$.
 For the special case
when all $f_k = D_{N_k}$, denote the matrices $M_k'$ by $M_k$, and denote
$\psi$ by $\psi_0$. It is clear
that $M_k$ has entries in $\{0,1\}$, $M_k'$ has entries in $\{-1,0, 1\}$,
$M_k'=M_k$ mod 2, and that in each column of $M_k$ the $1$ entries occur in an interval.
Hence $\det T(n,\psi_0) \leq \det T(n,\psi)$, by Corollary 2.3. Taking $n$th
roots and letting $n \to \infty$
on both sides of the latter inequality completes the proof, by the above theorem
of Szeg\"o.
\endproof
%

\centerline{\bf \S 3. Remarks and Questions. }
\medskip
{\bf 3.1.} The proof of (\ref{det3}) in Corollary 2.3 actually proceeded by
establishing the formally stronger coefficient-wise inequality
\begin{equation}
\label{det5}
D_n(A_1, A_2, \dots , A_n) \leq D_n(A_1', A_2', \dots , A_n').
\end{equation}
(This may be called a ``polarization" of its special case $\det A_1 \leq \det A_1'$,
whence the title of this paper.)
We want to note now that, as in (\ref{log}) of Theorem 1.1, the stronger inequality (\ref{det5}) may
also be written in terms of integrals
over circles, that is, in the case of Toeplitz matrices $A_k$ and  $A_k'$
generated as before by the $M_k$ of some $D_{N_k}(e^{i\theta})$ and the $M_k'$ of the
$f_k(e^{i\theta})$ respectively, for $k=1, \dots n$. In such a case one can invoke the
polarized version of the identity of Heine and Szeg\"o \cite[Theorem 1]{bum},
by which the right hand side of (\ref{det5}) becomes the multiple integral
\begin{equation}
\label{det6}
\frac{1}{(2\pi)^n}\int_0^{2\pi} \dots \int_0^{2\pi} \left (\prod_{k=1}^n |f_{k}(e^{i\theta_k})|^2 \right)
\ \Delta(\theta_1, \dots, \theta_n) d\theta_1 \dots d\theta_n
\end{equation}
where $\Delta(\theta_1, \dots, \theta_n)=
{\displaystyle \prod_{1 \leq p < q \leq n} }|e^{i\theta_p} - e^{i\theta_q}|^2$.
(The left hand side of (\ref{det5}) is of course the same integral with the $f_k$
replaced by the $D_{N_k}$.)

We note that a similar kind of polarization occurred in
the original motivating work of Hardy, Littlewood, and Gabriel
(\cite{HL28}, \cite{G}) on $L_p$-norm results concerning $\mathcal{ \widetilde{L}}$.
There the authors proved a rearrangement theorem for all polynomials,
with arbitrary coefficients, which can be specialized to our
context as the following result: For even integers $p=2s \geq 2$ one has
 $||f(e^{i\theta})||_p \leq ||D_N(e^{i\theta})||_p $,
 for all $f \in \mathcal{ \widetilde{L}}(N)$ having complex coefficients of modulus $1$ or $0$.
  Their proof proceeded via a stronger, polarized version, essentially that
\begin{equation}
\label{pol}
\int_0^{2\pi} |f_1(e^{i\theta})|^2 |f_2(e^{i\theta})|^2 \dots |f_s(e^{i\theta})|^2 d\theta \leq
\int_0^{2\pi} |D_{N_1}(e^{i\theta})|^2 |D_{N_2}(e^{i\theta})|^2 \dots |D_{N_s}(e^{i\theta})|^2 d\theta
\end{equation}
for any $N_k$ and $f_k \in \mathcal{ \widetilde{L}}(N_k)$ with coefficients of modulus $1$ or $0$.
By a stretch of the imagination, one could view these integrals as being similar
to those in (\ref{det6}) above. In fact, one simply replaces $\Delta$ by a singular
measure concentrated on the ``diagonal" $\theta_1= \dots = \theta_n$ of the $n$-torus.
It would thus be interesting to investigate what other ``weights'' in place of
$\Delta$ would yield true inequalities (in one direction or the other).

{\bf 3.2.} To what extent can the hypotheses that
$ f_k \in \mathcal{L}(N_k)$ be relaxed, in the inequalities proved in this paper?
For example, can we allow $f_k \in \mathcal{ \widetilde{L}}(N_k)$ ?
The strong term-wise inequality
(\ref{det4}) is trivially false in general for $f_k \in \mathcal{ \widetilde{L}}(N_k)$,
even when the polynomials $f_k$ have coefficients in $\{0,1\}$ (take $K=n=1, N_1=2$
and $f_1(z) = 1 + z^2$).
Going up one level of summation, does the inequality (\ref{det5}) hold whenever
$f_k \in \mathcal{ \widetilde{L}}(N_k)$ ? The answer is again no; numerical
work by the author has uncovered counterexamples to (\ref{det5}) with $n=7, \ N_1 = 2, \ f_1(z) =1+cz,
\ N_2=6=N_3 = \dots = N_7, \ f_2(z)
= 1+\sum_{j=1}^{5} c_jz^j =f_3(z)= \dots = f_7(z)$, with
some complex coefficients satisfying $|c_j|=1=|c|$. 
What if only $\pm 1$ coefficients
(and gaps) are allowed in the definition of $\mathcal{\widetilde{L}}(N)$, that is
$f_k(z) = \sum_{j=1}^{N_k} \pm z^{n_j}$ ?

{\bf 3.3.} We finally mention questions on the subject of sharp
$L_p$ inequalities for  real $p$:
It seems natural to ask whether (\ref{L24})
holds for all $p \geq 2$, and for all $f_k \in \mathcal{ \widetilde{L}}(N_k)$
with coefficients of modulus $1$ or $0$, or at least for all $f_k \in \mathcal{{L}}(N_k)$,
and whether matrix versions hold, such as (\ref{p2}) for all $p \geq 1$.
In connection with the integer $p$ case, a natural matrix version of (\ref{pol}) would be
to ask whether
\begin{equation}
\label{pol2}
{\rm trace}(A_1'A_2' \dots A_s') \leq {\rm trace}(A_1A_2 \dots A_s)
\end{equation}
for the matrices discussed above in 3.1, when $f_k \in \mathcal{ \widetilde{L}}(N_k)$
have coefficients of modulus $1$ or $0$. By taking
absolute values of everything inside the trace, one sees that (\ref{pol2}) 
holds trivially for $f_k \in \mathcal{ {L}}(N_k)$, and one could envision
adapting the rearrangement arguments of \cite{G} to prove it
for $\mathcal{ \widetilde{L}}$. However, this would not imply anything for
non-integer $p$'s. Instead, to prove (\ref{p2}) for non-integers $p >2$,
one may need to study some further homogeneous polynomials in the matrix entries,
generalizing elementary symmetric polynomials (such as in \cite{k3} and \cite{k4}).




\bigskip

\end{document}